\documentclass[article,12pt]{elsarticle}

\usepackage[small,bf,hang]{caption} 

\usepackage{graphicx}
\usepackage{subfig}

\usepackage{hyperref}

\journal{Discrete Applied Mathematics}

\newtheorem{theorem}{Theorem}[section]

\newtheorem{conjecture}[theorem]{Conjecture}
\newtheorem{rconjecture}[theorem]{Refuted Conjecture}
\newtheorem{observation}[theorem]{Observation}

\graphicspath{{figures/}}

\begin{document}

\begin{frontmatter}

\title{A counterexample to the pseudo 2-factor \\isomorphic graph conjecture}

\author[gent]{Jan Goedgebeur\fnref{fwo}}
\ead{Jan.Goedgebeur@UGent.be}

\address[gent]{Department of Applied Mathematics, Computer Science \& Statistics\\
  Ghent University\\
Krijgslaan 281-S9, \\9000 Ghent, Belgium\\ }

\fntext[fwo]{Supported by a Postdoctoral Fellowship of the Research Foundation Flanders (FWO).}

\begin{abstract}
A graph $G$ is \textit{pseudo 2-factor isomorphic} if the parity of the number of cycles in a 2-factor is the same for all 2-factors of $G$. Abreu et al.~\cite{abreu_08} conjectured that $K_{3,3}$,  the Heawood graph and the Pappus graph are the only essentially 4-edge-connected pseudo 2-factor isomorphic cubic bipartite graphs (Abreu et al., Journal of Combinatorial Theory, Series B, 2008, Conjecture 3.6).

Using a computer search we show that this conjecture is false by constructing a counterexample with 30 vertices. We also show that this is the only counterexample up to at least 40 vertices.

A graph $G$ is \textit{2-factor hamiltonian} if all 2-factors of $G$ are hamiltonian cycles. Funk et al.~\cite{funk_03} conjectured that every 2-factor hamiltonian cubic bipartite graph can be obtained from $K_{3,3}$ and the Heawood graph by applying repeated star products (Funk et al., Journal of Combinatorial Theory, Series B, 2003, Conjecture 3.2). We verify that this conjecture holds up to at least 40 vertices.
\end{abstract}


\end{frontmatter}

\section{Introduction and preliminaries}

All graphs considered in this paper are simple and undirected. Let $G$ be such a graph. We denote the vertex set of $G$ by $V(G)$ and the edge set by $E(G)$.

\begin{figure}[h!t]
    \centering
    \subfloat[]{\label{fig:heawood}\includegraphics[width=0.35\textwidth]{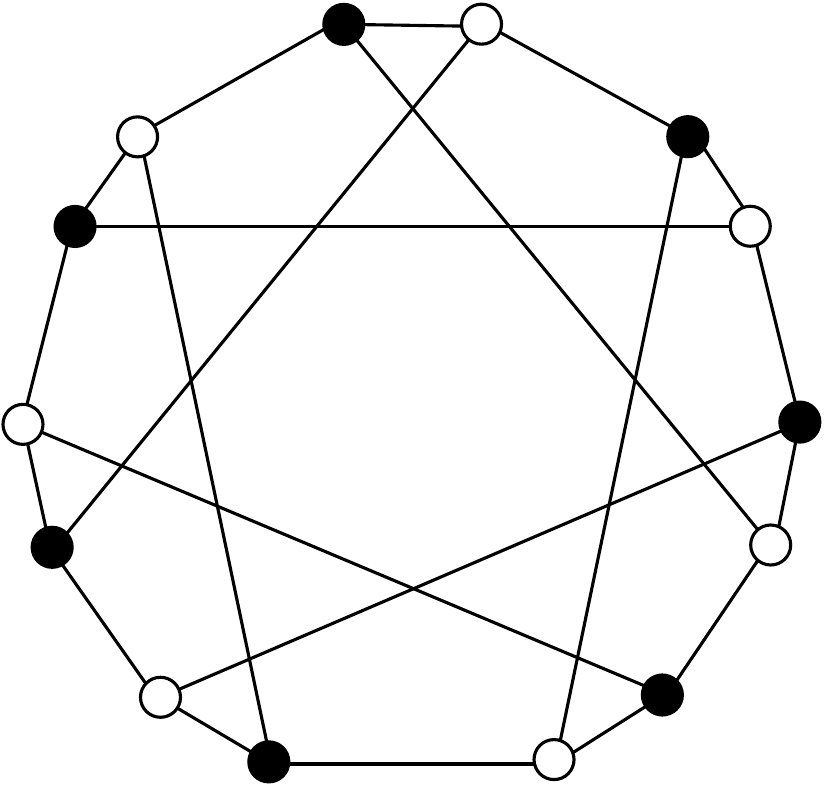}}\qquad\ \qquad
    \subfloat[]{\label{fig:pappus}\includegraphics[width=0.35\textwidth]{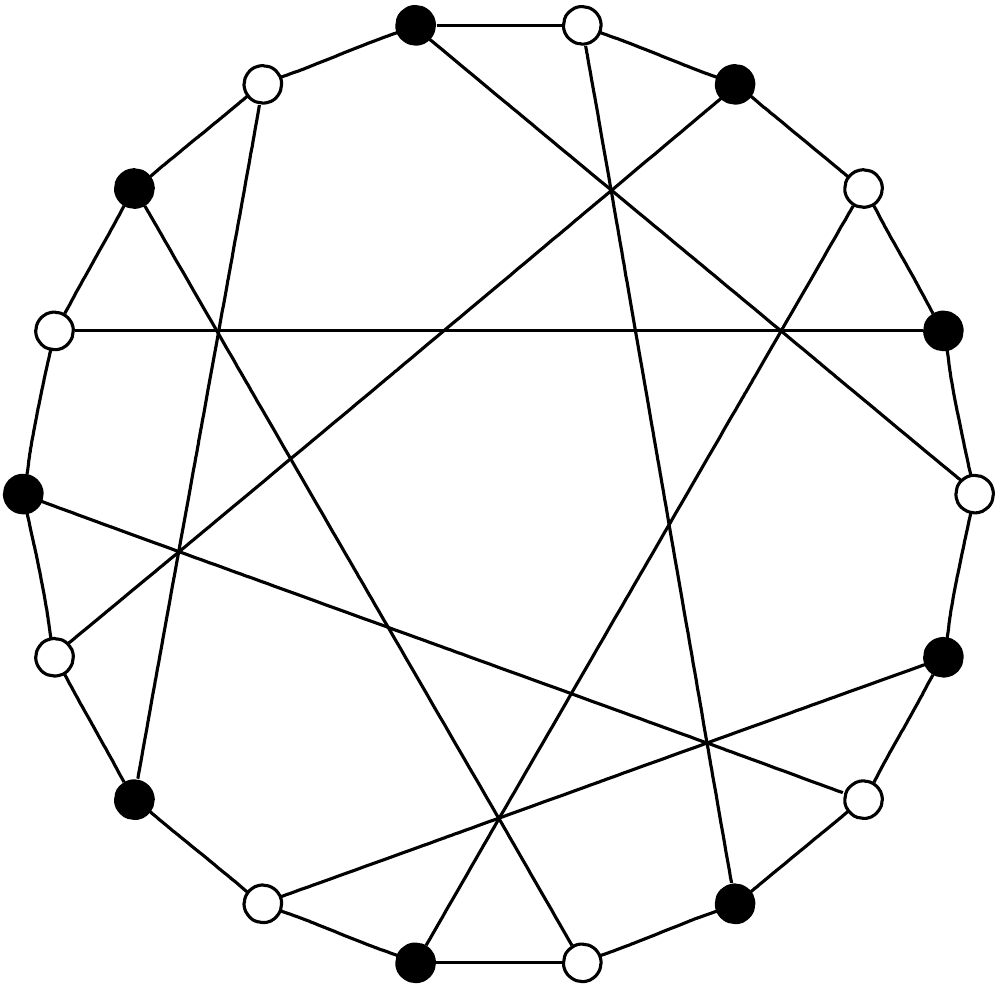}}
    \caption{The Heawood graph (a) and the Pappus graph (b).}
    \label{fig:cubic_graphs}
\end{figure}

A graph $G$ is \textit{2-factor hamiltonian} if all 2-factors of $G$ are hamiltonian cycles. Examples of 2-factor hamiltonian graphs include $K_4, K_5, K_{3,3}$ and the Heawood graph (see Figure~\ref{fig:heawood}). Funk et al.~\cite{funk_03} have shown that 2-factor hamiltonian $k$-regular bipartite graphs only exist when $k \in \{2,3\}$. They also constructed an infinite family of cubic bipartite 2-factor hamiltonian graphs obtained by applying repeated \textit{star products} to $K_{3,3}$ and the Heawood graph.

Given two cubic graphs $G_1, G_2$. A graph $G$ is a \textit{star product} of $G_1$ and $G_2$ if and only if there is an $x \in V(G_1)$ with neighbours $x_1,x_2,x_3$ in $G_1$ and an $y \in V(G_2)$ with neighbours $y_1,y_2,y_3$ in $G_2$ such that $G = (G_1 - x) \cup (G_2 - y) \cup \{(x_1,y_1),(x_2,y_2),(x_3,y_3)\}$.

Funk et al.\ conjectured that every 2-factor hamiltonian cubic bipartite graph belongs to their their infinite family of 2-factor hamiltonian cubic bipartite graphs, i.e.:

\begin{conjecture}[Funk et al., Conjecture 3.2 in~\cite{funk_03}]
Let $G$ be a 2-factor hamiltonian $k$-regular bipartite graph. Then either $k = 2$ and $G$ is a circuit or $k = 3$ and $G$ can be obtained from $K_{3,3}$ and the Heawood graph by repeated star products.
\label{2factorhamiltonian}
\end{conjecture}

As already mentioned in~\cite{funk_03}, it follows from~\cite{labbate_01} that a smallest counterexample to this conjecture is cubic and cyclically 4-edge-connected and from~\cite{labbate_02} that it has girth at least 6. So to prove Conjecture~\ref{2factorhamiltonian}, it would be sufficient to prove the following conjecture:

\begin{conjecture}[Funk et al.~\cite{funk_03}]
The Heawood graph is the only 2-factor hamiltonian cyclically 4-edge-connected cubic bipartite graph of girth at least~6.
\label{2factorhamiltonian2}
\end{conjecture}

Abreu et al.~\cite{abreu_08} extended these results on 2-factor hamiltonian graphs to the more general family of pseudo 2-factor isomorphic graphs. A graph $G$ is \textit{pseudo 2-factor isomorphic} if the parity of the number of cycles in a 2-factor is the same for all 2-factors of $G$. 

Clearly all 2-factor hamiltonian graphs are also pseudo 2-factor isomorphic. An example of a pseudo 2-factor isomorphic graph which is not 2-factor hamiltonian is the Pappus graph (see Figure~\ref{fig:pappus}), as the cycle sizes of its 2-factors are (6,6,6) and (18) (see~\cite{abreu_08}).

Abreu et al.~\cite{abreu_08} proved that pseudo 2-factor isomorphic $k$-regular bipartite graphs only exist when $k \in \{2,3\}$. They also constructed an infinite family of pseudo 2-factor isomorphic cubic bipartite graphs based on $K_{3,3}$, the Heawood graph and the Pappus graph and conjectured that these are the only pseudo 2-factor isomorphic cubic bipartite graphs:

\begin{rconjecture}[Abreu et al., Conjecture 3.5 in~\cite{abreu_08}]
Let $G$ be a 3-edge-connected cubic bipartite graph. Then $G$ is pseudo 2-factor isomorphic if and only if $G$ can be obtained from $K_{3,3}$, the Heawood graph or the Pappus graph by repeated star products.
\label{pseudo2factorisomorphic}
\end{rconjecture}

Note that the pseudo 2-factor isomorphic graphs obtained by such a star product have a non-trivial 3-edge-cut. An edge cut $E_0$ of a graph $G$ is \textit{essential} if $G - E_0$ has at least two non-trivial components.
A graph $G$ is \textit{essentially 4-edge-connected} if $G$ does not have an essential edge cut $E_0$ with $|E_0| < 4$. Therefore Conjecture~\ref{pseudo2factorisomorphic} can only hold if the following conjecture also holds:

\begin{rconjecture}[Abreu et al., Conjecture 3.6 in~\cite{abreu_08}]
Let $G$ be an essentially 4-edge-connected pseudo 2-factor isomorphic cubic bipartite graph. Then $G$ must be $K_{3,3}$, the Heawood graph or the Pappus graph.
\label{pseudo2factorisomorphic2}
\end{rconjecture}

It follows from Theorem~3.10 in~\cite{abreu_08} that $K_{3,3}$ is the only essentially 4-edge-connected pseudo 2-factor isomorphic cubic bipartite graph with girth 4, so a counterexample to Conjecture~\ref{pseudo2factorisomorphic2} must have girth at least 6.

Abreu et al.\ partially proved this conjecture for \textit{irreducible Levi graphs} (see~\cite{abreu_12} for details).

In the next section we describe the results of a computer search for cubic bipartite graphs of girth at least 6. This allowed us to verify Conjecture~\ref{2factorhamiltonian2} up to 40 vertices. It also yielded one counterexample with 30 vertices to Conjecture~\ref{pseudo2factorisomorphic2}. This is the only counterexample up to at least 40 vertices.

\section{Testing and results}
Using the program \textit{minibaum}~\cite{brinkmann_96} we generated all cubic bipartite graphs with girth at least 6 up to 40 vertices and all cubic bipartite graphs with girth at least 8 up to 48 vertices. The counts of these graphs can be found in Table~\ref{table:cubic_bip_graphs}. Some of these graphs can be downloaded from \url{http://hog.grinvin.org/Cubic}

\begin{table}[ht!]
\centering
\begin{tabular}{| c | r | r |}
\hline 
Vertices & Girth at least 6 & Girth at least 8\\
\hline 
14  &  1  &  \\
16  &  1  &  \\
18  &  3  &  \\
20  &  10  &  \\
22  &  28  &  \\
24  &  162  &  \\
26  &  1 201  &  \\
28  &  11 415  &  \\
30  &  125 571  &  1\\
32  &  1 514 489  &  0\\
34  &  19 503 476  &  1\\
36  &  265 448 847  &  3\\
38  &  3 799 509 760  &  10\\
40  &  57 039 155 060  &  101\\
42  &  ?  &  2 510\\
44  &  ?  &  79 605\\
46  &  ?  &  2 607 595\\
48  &  ?  &  81 716 416\\
\hline
\end{tabular}

\caption{Counts of cubic bipartite graphs with girth at least 6 or girth at least~8.}

\label{table:cubic_bip_graphs}
\end{table}

We then implemented a program which tests if a given graph is pseudo 2-factor isomorphic and applied it to the generated cubic bipartite graphs. This yielded the following results:

\begin{observation}
There is exactly one essentially 4-edge-connected pseudo 2-factor isomorphic graph different from the Heawood graph and the Pappus graph among the cubic bipartite graphs with girth at least 6 with at most 40 vertices.
\end{observation}

\begin{observation}
There is no essentially 4-edge-connected pseudo 2-factor isomorphic graph among the cubic bipartite graphs with girth at least 8 with at most 48 vertices.
\end{observation}

\begin{figure}[h!t]
	\centering
	\includegraphics[width=0.5\textwidth]{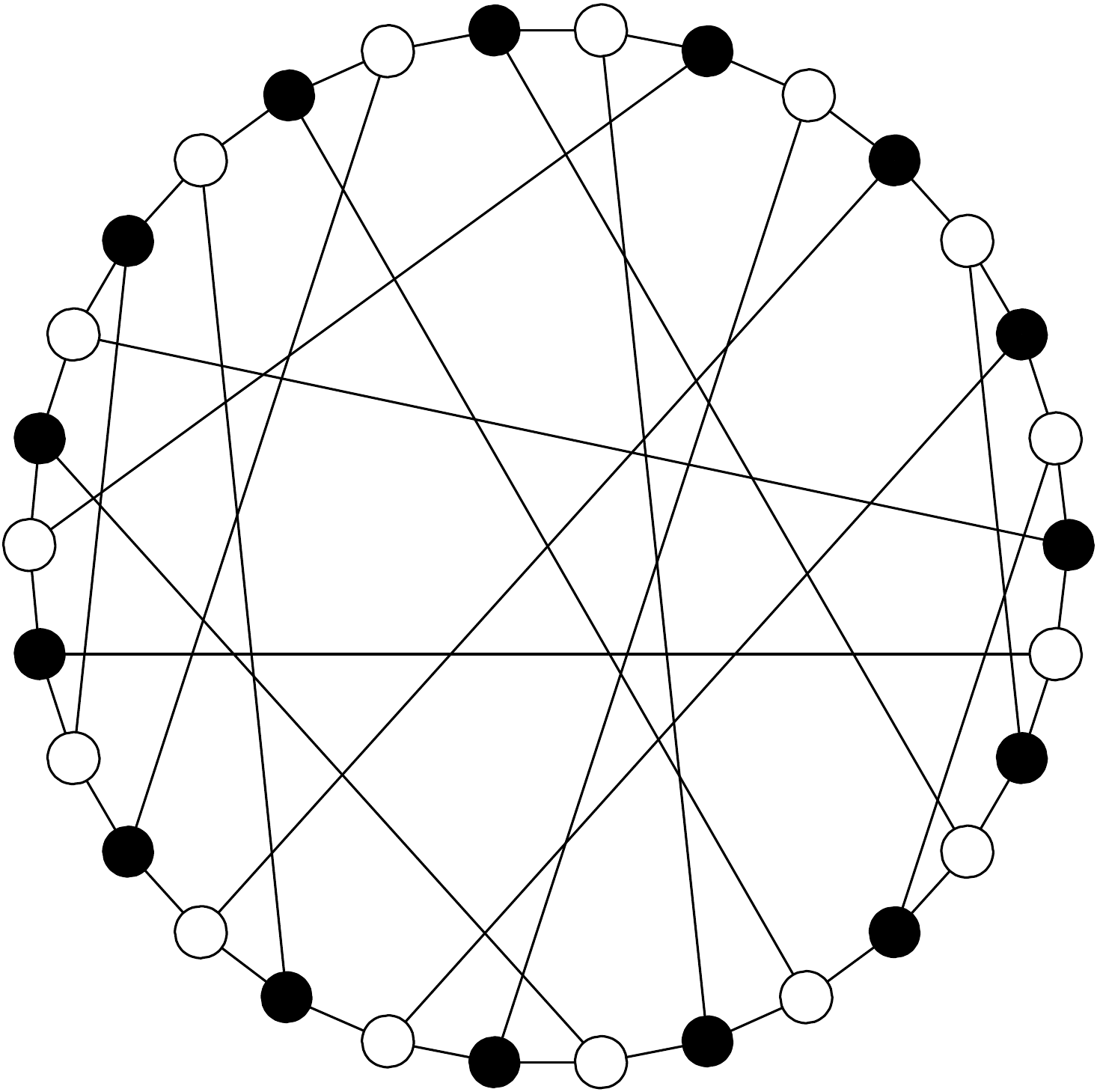}
	\caption{The pseudo 2-factor isomorphic cubic bipartite graph $\mathcal{G}$ with 30 vertices which is a counterexample to Conjecture~\ref{pseudo2factorisomorphic2}.}
	\label{fig:counterexample}
\end{figure}

This implies that Conjecture~\ref{pseudo2factorisomorphic2} (and consequently also Conjecture~\ref{pseudo2factorisomorphic}) is false. The counterexample has 30 vertices and there are no additional counterexamples up to at least 40 vertices and also no counterexamples among the cubic bipartite graphs with girth at least 8 up to at least 48 vertices. The counterexample (which we will denote by $\mathcal{G}$) is shown in Figure~\ref{fig:counterexample} and its adjacency list can be found in Table~\ref{table:counterexample}.
$\mathcal{G}$ can also be obtained from the \textit{House of Graphs}~\cite{hog} by searching for the keywords ``pseudo 2-factor isomorphic * counterexample'' where it can be downloaded and several of its invariants can be inspected.

$\mathcal{G}$ has cyclic edge-connectivity 6, automorphism group size 144, is not vertex-transitive, has 312 2-factors and the cycle sizes of its 2-factors are: (6,6,18), (6,10,14), (10,10,10) and (30). The list of all perfect matchings of this graph (and the corresponding 2-factors) can be found online in~\cite{pseudo2factor-site}.

Since all 2-factor hamiltonian graphs are pseudo 2-factor isomorphic and $\mathcal{G}$ is not 2-factor hamiltonian, this implies the following:

\begin{observation}
Conjecture~\ref{2factorhamiltonian2} holds up to at least 40 vertices and holds for cubic bipartite graphs with girth at least 8 up to at least 48 vertices.
\end{observation}

\begin{table}[ht!]
\centering
\begin{tabular}{r r r r}
0: & 1 & 13 & 29 \\
1: & 0 & 2 & 26 \\
2: & 1 & 3 & 21 \\
3: & 2 & 4 & 28 \\
4: & 3 & 5 & 19 \\
5: & 4 & 6 & 22 \\
6: & 5 & 7 & 15 \\
7: & 6 & 8 & 24 \\
8: & 7 & 9 & 27 \\
9: & 8 & 10 & 18 \\
10: & 9 & 11 & 25 \\
11: & 10 & 12 & 20 \\
12: & 11 & 13 & 17 \\
13: & 0 & 12 & 14 \\
14: & 13 & 15 & 23 \\
15: & 6 & 14 & 16 \\
16: & 15 & 17 & 29 \\
17: & 12 & 16 & 18 \\
18: & 9 & 17 & 19 \\
19: & 4 & 18 & 20 \\
20: & 11 & 19 & 21 \\
21: & 2 & 20 & 22 \\
22: & 5 & 21 & 23 \\
23: & 14 & 22 & 24 \\
24: & 7 & 23 & 25 \\
25: & 10 & 24 & 26 \\
26: & 1 & 25 & 27 \\
27: & 8 & 26 & 28 \\
28: & 3 & 27 & 29 \\
29: & 0 & 16 & 28 \\
\end{tabular}

\caption{The adjacency list of the graph $\mathcal{G}$ which is a counterexample to Conjecture~\ref{pseudo2factorisomorphic2}.}

\label{table:counterexample}
\end{table}

Independent checks are very important for computational results to minimize the chance of errors.

The correctness of the program \textit{minibaum}~\cite{brinkmann_96} which was used to generate cubic bipartite graphs has been verified before, e.g.\ by comparing it to other generators for cubic graphs such as~\cite{brinkmann_01} or~\cite{meringer_99}.

We also implemented two independent algorithms to test if a given graph is pseudo 2-factor isomorphic. One implementation searches for 2-factors by building disjoint cycles and the other constructs perfect matchings. One of the implementations is based on a program which has already been extensively used and tested in~\cite{snark-paper}.

We applied both programs to the lists of generated cubic bipartite graphs and verified that they yielded the same results. We also modified the two programs to count the 2-factors in a given graph. We applied these modified programs to several graphs and verified that both programs found the same number of 2-factors. 

Since all results of the two independent algorithms are in complete agreement, we believe that this is strong evidence for the correctness of our implementations and results.

\section*{Acknowledgements}
Most computations for this work were carried out using the Stevin Supercomputer Infrastructure at Ghent University.


\bibliographystyle{plain}
\bibliography{references}

\end{document}